# Probabilistic approach to the distribution of primes and to the proof of Legendre and Elliott-Halberstam conjectures

VICTOR VOLFSON

ABSTRACT. Probabilistic models for the distribution of primes in the natural numbers are constructed in the article. The author found and proved the probabilistic estimates of the deviation $R(x) = |\pi(x) - Li(x)|$. The author has analyzed the probabilistic models of the distribution of primes in the natural numbers and affirmed the validity of the probabilistic estimates of proved deviations $R(x)$ stronger than the estimates made under the assumption of Riemann conjecture. Legendre's conjecture was proved in this paper with probability arbitrarily close to 1 based on the probability estimates. Probabilistic models for the distribution of primes in the arithmetic progression $ki + l$, $(k,l) = 1$ are also built in this paper. The author has proved the probability estimates for the deviation $R(x,k,l) = |\pi(x,k,l) - Li(x)/\varphi(k)|$. The author has analyzed the probability models of the distribution of primes in the arithmetic progression and affirmed the validity of probabilistic estimates of proved deviations $R(x,k,l)$ stronger than the estimates made under the assumption of the extended Riemann conjecture. Elliott-Halberstam conjecture $\sum_{1 \leq k \leq x^a} \{max_{(k,l)=1} R(x,k,l)\} \leq Cx/\ln^A(x)$ was proved in this paper with probability arbitrarily close to 1 for all $0 < a < 1$ based on the probability estimates.

## 1. INTRODUCTION

The probabilistic methods have recently been used to model the distribution of prime numbers. Hardy, Littlewoods and Kramer were the first to make it in their conjectures. These conjectures showed that when modeling the distribution of primes probabilistic methods provide stronger results than the theory of functions of a complex variable and therefore cannot be underestimated.

The task of studying the distribution of values of real arithmetic functions can be associated with the theory of summation of random variables in some cases [1]. It turns out that the question of the asymptotic behavior of these functions can be reduced to the known limit theorems of the probability theory.

Keywords: probabilistic model, limit distribution of the random variable, normal distribution, probability estimates, prime numbers, Riemann conjecture, Legendre conjecture, Cramer conjecture, Elliott-Halberstam conjecture.

It is interesting to see the probabilistic assessment of real arithmetic functions of the number of primes in the natural numbers and arithmetic progression. This question is discussed in this paper.

This article considers several probability spaces such that the actual arithmetic function of the number of primes in the natural numbers and arithmetic progression are relatively simple sums of random variables.

The asymptotic law of prime numbers proved for the first time using the theory of functions of a complex variable is well known in the following form:

$$\pi(x) \sim x/\ln(x), \tag{1.1}$$

where $\pi(x)$ is the number of primes not exceeding a positive integer $x$. The analytical approach of Chebyshev was an early attempt to estimate the value:

$$R(x) = |\pi(x) - x/\ln(x)|. \tag{1.2}$$

In [2] the estimate:

$$ax/\ln(x) < \pi(x) < bx/\ln(x) \tag{1.3}$$

was shown to be valid. Chebyshev showed that $a = 0.921$ and $b = 1.106$. Values more close to 1 were obtained in subsequent papers. From (1.3) it follows:

$$(a-1)x/\ln(x) < \pi(x) - x/\ln(x) < (b-1)x/\ln(x). \tag{1.4}$$

However, estimates of type (1.4) are inaccurate:

$$R(x) = O(x/\ln(x)). \tag{1.5}$$

The asymptotic law of the prime numbers gives more accurate value $\pi(x)$ using the integral logarithm:

$$\pi(x) \sim Li(x). \tag{1.6}$$

The accuracy of the formula (1.6) is given by:

$$|\pi(x) - Li(x)| = O(x^{1/2} \cdot \ln(x)), \tag{1.7}$$

if the Riemann hypothesis is true.

There is also Legendre's formula:

$$\pi(x) \sim x/\ln(x) + B. \tag{1.8}$$

But the formula (1.8) is less accurate than equation (1.6). This raises the question whether it possible to get a better probabilistic assessment of the accuracy of the formula:

$$\pi(x) = x/\ln(x) + x \cdot o(1/ln(x)). \tag{1.9}$$

To answer this question, we consider probabilistic models of the distribution of primes in the natural numbers.

## 2. PROBABILISTIC MODELS OF THE DISRIBUTION OF PRIMES IN THE NATURAL NUMBERS

Consider the first probabilistic model. Let there be $x$ balls indistinguishable to the touch. Number them by consecutive natural numbers from 1 to $x$ and put them all in a box.

Select one ball from the box at random. If its number belongs to a pre-selected integer positive injective sequence, then we consider the event as "success", and if the number does not belong to the selected sequence, then we consider this event as "failure".

We assume that the probability of the successful event is $p$. Accordingly, the probability of failure events will be equal to $1-p$. We introduce the random variable $I_1$ – an indicator of the success of the event. The value of the random variable is equal to $I_1 = 1$, if it is success and value is equal to $I_1 = 0$, if it is failure. Return the first ball in the box, toss balls in the box and choose the second ball from the box at random. If the ball number belongs to the selected sequence, then assign to the random variable the value $I_2 = 1$. If the ball number does not belong to the selected sequence, then assign to the random variable the value $I_2 = 0$. Then return the second ball in the box, too. Repeat it $x$ times. Since the balls are each time returned to the box and every next ball is resampled in the same conditions, the random variables $I_i$ are independent. Thus, we obtain a sequence of independent random variables – indicators of success of the event: $I_1, I_2, ... I_x$.

The expectation of a random variable $I_i$ is equal to:

$$M(I_i) = p \cdot 1 + (1-p) \cdot 0 = p. \tag{2.1}$$

The variance of a random variable $I_i$ is equal to:

$$D(I_i) = (1-p)^2 \cdot p + p^2(1-p) = p(1-p). \tag{2.2}$$

Let us consider the random variable equal to the sum:

$$I(x) = \sum_{i=1}^{x} I_i. \qquad (2.3)$$

Define the characteristics of the random variable $I(x)$. The expectation of the random variable $I(x)$ based on (2.1) and (2.2) and the linearity of expectations is equal to:

$$M(I(X)) = M(\sum_{i=1}^{x} I_i) = \sum_{i=1}^{x} M(I_i) = x \cdot p. \qquad (2.4)$$

Based on (2.2), (2.3) and independency of the random variables $I_i$, the variance of the random variable $I(x)$ is equal to;

$$D(I(x)) = D(\sum_{i=1}^{x} I_i) = \sum_{i=1}^{x} D(I_i) = x \cdot p(1-p). \qquad (2.5)$$

Thus, we have the independent, identically distributed random variables $I_1, I_2, ... I_x$ with bounded dispersion. Therefore, the random variable $I(x) = \sum_{i=1}^{x} I_i$ has a binomial distribution. The limiting distribution for the random variable $I(x)$ for the fixed probability $p$ is the normal distribution. Based on Moivre-Laplace theorem [3], the following relation is held:

$$\lim_{x \to \infty} \{P(|I(x) - M(I(x))| < C\sqrt{D(I(x))})\} = F(C), \qquad (2.6)$$

where $P()$ the probability of the event is specified in parentheses and $F(C)$ is the function of the module of the standard normal distribution at point C.

Let us substitute the characteristics of the random variable $I(x)$ from (2.4) and (2.5) into the expression (2.6) and obtain:

$$\lim_{x \to \infty} \{P(|I(x) - x \cdot p| < C\sqrt{x \cdot p(1-p)})\} = F(C). \qquad (2.7)$$

The value $F(C)$ tends rapidly to 1 with increasing $C$, based on the properties of the normal distribution. Thus, we can choose such value $C$ that the probability of the event:

$$|I(x) - x \cdot p| < C\sqrt{x \cdot p(1-p)} \qquad (2.8)$$

is arbitrarily close to 1.

This is an advantage of probabilistic assessments as it gives the opportunity to select a particular value $C$. For example, formula (2.8) has the benefit, as compared to the formula (1.7), in which the accuracy cannot be evaluated.

Denote by $\pi(f, A, B)$ the number of members of the integer, positive, injective sequence $f(n)$ on the interval $[A, B)$. It was proved in [4] that the density of the sequence $f(n)$ of the natural numbers within the interval $[A, B)$ defined by the formula:

$$\pi(f, A, B) / (B - A) \qquad (2.9)$$

is a finite probability measure.

The sequence of prime numbers is an integer, positive, injective, i.e. satisfies the above conditions.

On the basis of the asymptotic law of prime numbers and (2.9), it was also shown in [4] that the probability of choosing a prime within the interval of the natural numbers $[2, x)$ is equal to:

$$p = 1/\ln(x) + o(1/\ln(x)). \qquad (2.10)$$

The random variable $I(x)$, in this case, can be considered as the number of primes not exceeding a natural number $x$ - $\pi(x)$ [1]. Thus the formula (2.7) for the sequence of primes and large x may be written on the basis of (2.10):

$$P(|\pi(x) - x \cdot (1/\ln(x) + o(1/\ln(x)))| <$$

$$< C\sqrt{x \cdot (1/\ln(x) + o(1/\ln(x)))(1 - 1/\ln(x) - o(1/\ln(x)))}) \approx F(C). \qquad (2.11)$$

On the basis of (2.11), for large values $x$ we can choose value $C$ such that the probability of an event:

$$|\pi(x) - x \cdot (1/\ln(x) + o(1/\ln(x)))| < C\sqrt{x \cdot (1/\ln(x) + o(1/\ln(x)))(1 - 1/\ln(x) - o(1/\ln(x)))} \qquad (2.12)$$

is arbitrarily close to 1.

Let us analyze formula (2.12). It gives the upper bound $\pi(x)$ for the deviation of the number of primes not exceeding $x$ from the value $x/\ln(x) + x \cdot o(1/\ln(x))$, i.e. the estimate of accuracy for the formula (1.9).

Consider this probabilistic model in the case when $o(1/\ln(x))$ is the function:

$$f(x) = \sum_{i=1}^{\infty} (i-1)/\ln^i(x) = Li(x)/x - 1/\ln(x). \quad (2.13)$$

Substitute function (2.13) into the formula for the variance of the random variable $I(x)$ in the first probabilistic model and obtain:

$$D_1(I(x)) = Li(x)(1 - Li(x)/x) = Li(x) - Li^2(x)/x = \int_2^x \frac{dt}{\ln(t)} - \frac{(\int_2^x \frac{dt}{\ln(t)})^2}{x}. \quad (2.14)$$

By substituting (2.14) into (2.12) we find that the probability of the event:

$$|\pi(x) - Li(x)| < C\sqrt{Li(x) - Li^2 x/x} \quad (2.15)$$

is arbitrarily close to 1.

Let us analyze this probabilistic model. After sampling (in this model) the ball returns to the box again. Therefore, in this model it is possible to choose the same ball several times. This does not happen when we calculate the number of members of the sequence in the range of the natural numbers from 1 to $x$ in the real situation.

Consider the second probabilistic model of the distribution of primes, which is free from this drawback. It will be based on the of Cramer's model.

Let $(P_n)$ be an infinite series of urns containing black and white balls? The chance of drawing a white ball from $U_n$ being $p_i$, while the composition of $U_1$ and $U_2$ may be arbitrary chosen. We now assume that one ball is drawn from each urn, so that an infinite series of alternately black and white balls is obtained. If $P_n$ denotes the number of the urn from which the $n$-th white ball in the series is drawn? The number $P_1, P_2, ....$ will from increasing sequences of integers, And we shall consider the class $C$ of all possible sequences $x$. Obviously the sequence $S$ of ordinary prime numbers $(p_n)$ belongs to this class.

We shall denote by $I(x)$ the number of those $P_i$, which are $\leq x$. We denote by $I_i$ the random variable that takes the value 1 if a white ball gets from the $i$-th urn and the value 0 otherwise. Thus we have $I(x) = \sum_{i=1}^{x} I_i$.

Thus, the probability model is free from the above drawbacks.

We shall find the characteristics of random variables of the second probabilistic model. The expectation of the random variable $I_i$ is equal to:

$$M(I_i) = 1 \cdot p_i + 0 \cdot (1-p_i) = p_i. \qquad (2.16)$$

The variance of the random variable $I_i$ is equal to:

$$D(I_i) = (1-p_i)^2 \cdot p_i + (p_i)^2 \cdot (1-p_i) = p_i(1-p_i). \qquad (2.17)$$

Define characteristics $I(x)$. Based on (2.16) the expectation $I(x)$ is equal to:

$$M(I(X)) = M(\sum_{i=1}^{x} I_i) = \sum_{i=1}^{x} M(I_i) = \sum_{i=1}^{x} p_i. \qquad (2.18)$$

The variance of the random variable $I(x)$ (based on independent random variables $I_i$ and formulas (2.17) is equal to:

$$D(I(x)) = D(\sum_{i=1}^{x} I_i) = \sum_{i=1}^{x} D(I_i) = \sum_{i=1}^{x} p_i(1-p_i) = \sum_{i=1}^{x} p_i - \sum_{i=1}^{x} (p_i)^2. \qquad (2.19)$$

The limiting distribution of the random variable $I(x)$ is normal, based on Lyapunov theorem [3]. Therefore, based on formulas (2.6), (2.17) and (2.19), we obtain:

$$\lim_{x \to \infty} \{P(|I(x) - \sum_{i=1}^{x} p_i| < C\sqrt{\sum_{i=1}^{x} p_i - \sum_{i=1}^{x} (p_i)^2}\} = F(C), \qquad (2.20)$$

where $F(C)$ is the value of the module of the function of the standard normal distribution at point $C$. Formula (2.20) (for large values $x$) can be written as:

$$P(|I(x) - \sum_{i=1}^{x} p_i| < C\sqrt{\sum_{i=1}^{x} p_i - \sum_{i=1}^{x} (p_i)^2}) \approx F(C). \qquad (2.21)$$

The probability tends rapidly to 1 with increasing $x$, based on the properties of the normal distribution. Thus, we can choose such a value $C$, that probability (2.21) is arbitrarily close to 1.

We shall find the characteristics of the probabilistic model for the sequence of prime numbers. Let us make a heuristic assumption that the probability of a natural number $n$ to be a prime is equal to $1/\ln(n)$, then the probability of choosing the white ball in the Kramer's model is $p_i = 1/\ln(i)$ for $i > 2$. Based on Kramer random variable $I(x)$ can be regarded as the number of primes not exceeding a natural number $x$ - $\pi(x)$.

The expectation of the random variable $I(x)$ for the sequence of the prime numbers, based on formula (2.18) and the heuristic assumption, is equal to:

$$M(I(x)) = \sum_{i=2}^{x} 1/\ln(i) \approx \int_{2}^{x} \frac{dt}{\ln(t)}. \qquad (2.22)$$

The variance of the random variable $I(x)$ for the sequence of prime numbers, based on formula (2.19) and the assumption, is equal to:

$$D(I(x)) = \sum_{i=2}^{x} 1/\ln(i) - \sum_{i=2}^{x} (1/\ln(i))^2 \approx \int_{2}^{x} \frac{dt}{\ln(t)} - \int_{2}^{x} \frac{dt}{\ln^2(t)}. \qquad (2.23)$$

Based on (2.21) we obtain (for large $x$):

$$P(|I(x) - \int_{2}^{x} \frac{dt}{\ln(t)}| < C\sqrt{\int_{2}^{x} \frac{dt}{\ln(t)} - \int_{2}^{x} \frac{dt}{\ln^2(t)}}) \approx F(C), \qquad (2.24)$$

where $F(C)$ is the value of the function of the module of the standard normal distribution at point $C$. Thus, we can choose such value $C$, that probability (2.24) is arbitrarily close to 1.

The question arises concerning the approximation sign in formulas (2.22) and (2.23). There are assertions 1, 2, 3 for partial summation [5] based on the Euler-McLaurin formula.

Assertion 1. Assume that the function $F(x)$ is continuously differentiable and monotonically decreasing on the interval $[A, \infty)$ and $\lim_{x \to \infty} F(X) = 0$. Then

$$\sum_{i=0}^{n} F(A+i) - \int_{A}^{n} F(x)dx = C + O(F(n)), \text{ where } C = \lim_{n \to \infty}[\sum_{i=0}^{n} F(A+i) - \int_{A}^{n} F(x)dx].$$

Assertion 2. Suppose that there exists the function $F(x)$ on the interval $[A, \infty)$ with the following properties: 1. $\lim_{x \to \infty} F(x) = 0$. 2. It has derivatives of the desired order. 3. $F^{(2k-1)}(x) < 0$.

4. $|\frac{B_{2k}(2K)!F^{(2k+1)}(A)}{B_{2k+2}(2k+2)!F^{(2k-1)}(A)}| < 1$, where $B_n$ is $n$ - Bernoulli number. Then

$$C \leq F(A)/2 + |F'(A)|/12, \text{ where } C = \lim_{n \to \infty}[\sum_{i=0}^{n} F(A+i) - \int_{A}^{n} F(x)dx].$$

Assertion 3. The following estimate is satisfied for the function $F(x) = 1/\ln^k(x)$ on the interval $[A, \infty)$:

$$C < 0.6202 F(k+1), \text{ where } C = \lim_{n \to \infty}[\sum_{i=0}^{n} F(A+i) - \int_{A}^{n} F(x)dx].$$

We obtain the following estimate based on assertion 1:

$$\sum_{i=2}^{x} 1/\ln(i) - \int_{2}^{x} dt/\ln(t) = C_1 + O(1/\ln(x)).$$

$$\sum_{i=2}^{x} 1/\ln^2(i) - \int_{2}^{x} dt/\ln^2(t) = C_2 + O(1/\ln^2(x)).$$

From assertions 2 and 3 it follows that $C_1 < 0.8948$, $C_2 < 0.6783$.

These values are negligible for large $x$. The function $1/\ln(x) - 1/\ln^2(x)$ decreases monotonically only for values $x > 7$, and therefore assertion 1 is applicable to this interval. Let us calculate separately the sum and the integral on the interval from 2 to 7. The sum is equal to 0.117 and the integral to 0.7119. In view of the difference between the sum and the integral, overall estimate of the constant is equal to $|C_3| < 0.47$ for variance, i.e. it is also negligible.

3. ANALYSIS OF THE PROBABILISTIC MODELS OF THE DISTRIBUTION OF THE PRIME NUMBERS

Compare the variance of the first probabilistic model defined by formula (2.14) with the variance of the second probabilistic model defined by formula (2.23). We use a special case of Cauchy-Schwarz inequality for comparison:

$$(\int_{2}^{x} u(x)v(x)dx)^2 \leq \int_{2}^{x} u^2(x)dx \cdot \int_{2}^{x} v^2(x)dx. \tag{3.1}$$

The equation (3.1) is written as $\dfrac{(\int_{2}^{x} v(x)dx)^2}{x-2} \leq \int_{2}^{x} v^2(x)dx,$ (3.2)

if $u(x) = 1$ and $v(x) > 0$. We have:

$$\frac{(\int_{2}^{x} \frac{dt}{\ln(t)})^2}{x} \leq \int_{2}^{x} \frac{dt}{\ln^2(t)}, \tag{3.3}$$

based on (3.2) for large values $x$.

From (3.3) we obtain:

$$D_2(I(x)) \leq D_1(I(x)). \tag{3.4}$$

Then based on (3.4), we have the following upper bound:

$$D_2(I(x)) \leq D_1(I(x)) \leq Li(x), \tag{3.5}$$

where $Li(x) = \int_2^x \frac{dt}{\ln(t)}$.

The question naturally arises about the validity of these probability estimates as the probability estimates for deviations $r(x) = |\pi(x) - Li(x)|$ are stronger than estimate (1.7) made under the assumption of the Riemann conjecture.

To confirm these probability estimates we have to find already proved lower bound for the function $R(x)$ that has the order less than $\sqrt{D_2(I(x))} = \sqrt{\int_2^x \frac{dt}{\ln(t)} - \int_2^x \frac{dt}{\ln^2(t)}}$ starting from a certain value $x$.

The paper [6] contains such an assessment. Littlewoods showed that the value $R(x)$ cannot have the order less than:

$$\frac{\sqrt{x} \ln \ln \ln(x)}{\ln(x)} \tag{3.6}$$

as $x$ tends to infinity. Therefore, it suffices to compare the order of the square of the function (3.6) with the value $D_2(I(x)) = \int_2^x \frac{dt}{\ln(t)} - \int_2^x \frac{dt}{\ln^2(t)}$. Let us find the limit of the value:

$$A = \lim_{x \to \infty} \{ \frac{x \ln^2(\ln \ln(x))}{\ln^2(x)(\int_2^x \frac{dt}{\ln(t)} - \int_2^x \frac{dt}{\ln^2(t)})} \}. \tag{3.7}$$

We find that the value of (3.7) is equal to:

$$A = \lim_{x \to \infty} \{ \frac{\ln^2(\ln(x))}{\ln(x)} \} + 2\lim_{x \to \infty} \{ \frac{\ln \ln \ln(x)}{(\ln \ln(x)) \ln^2(x)} \} = 0, \tag{3.8}$$

(using L'Hospital's Rule and the derivative on the upper limit of integration), since both expressions in the sum tend to 0. Thus, the function (3.6) has the order of magnitude smaller than $\sqrt{D_2(I(x))} = \sqrt{\int_2^x \frac{dt}{\ln(t)} - \int_2^x \frac{dt}{\ln^2(t)}}$. It validates estimates made based on the first and second probabilistic models.

The second probability model is more accurate than the first probability model, as it does not have the disadvantages of the first one. On the other hand, based on (3.5), dispersion for the first probabilistic model (2.14) is greater than that for the second probabilistic model (2.23), so it can be used for the upper estimate in (2.15). Thus, the first probabilistic model is strictly proved (without assuming the second probabilistic model) and inequality (2.15) (for the first probabilistic model) is rigorously proved assertion.

Numerical calculations show that the inequality is at least satisfied for values $x < 10^9$. However, Littlewoods, as mentioned earlier, proved the theorem [6] and showed that this case is not always held. But this theorem does not find such a value $x$ that the inequality $\pi(x) > Li(x)$ is satisfied. Skewers have already answered this question. He found such value $x = \exp\exp\exp\exp(7.705)$ [7] and confirmed this fact.

The probabilistic models of the distribution of primes in the natural numbers listed above correspond to the specified Littlewoods theorem and prove that the inequality sign can be changed. It has been conclusively established that these deviations cancel each other [8]. This is consistent with the results obtained on the standard normal distribution of the deviation. Therefore, the term of "the distribution of prime numbers" is well founded in the terminology of the theory of probability.

## 4. APPLICATION OF THE PROBABILISTIC MODELS OF THE DISTRIBUTION OF PRIMES IN THE NATURAL NUMBERS

Suppose that the Riemann conjecture is true. Then (for a sufficiently large number $A > 0$ and $x$) there always exists at least one prime between the numbers [9]:

$$x \text{ and } x + Ax^{1/2}\ln(x). \tag{4.1}$$

Similarly, from (2.14) it follows that there always exists at least one prime between the values:

$$x \text{ and } x + C\sqrt{Li(x) - Li^2(x)/x} \tag{4.2}$$

with probability arbitrarily close to 1.

Based on (4.6), note that the condition:

$$\sqrt{Li(x) - Li^2(x)/x} < \sqrt{Li(x)} < \sqrt{x} . \tag{4.3}$$

is held. Thus, condition (4.2) is stronger than (4.1).

Recall that the Kramer conjecture assumes that there always exist primes between the numbers:

$$x \text{ and } x + B\ln^2(x). \tag{4.4}$$

for large values $x$ and appropriate $B > 0$. Therefore, from (4.4) we have:

$$x + Ax^{1/2}\ln(x) > x + C\sqrt{Li(x) - Li^2(x)/x} > x + B\ln^2(x) . \tag{4.5}$$

Thus, the condition (4.4) is stronger than (4.2), but (4.4) is only the conjecture.

Legendre conjecture, that there exists at least one prime between the squares of consecutive positive integers are well known.

Assertion 4. Legendre conjecture is held with probability arbitrarily close to 1.

Proof. It is known [9], that the Legendre conjecture will be proved, if we show, that there always exists at least one prime between some large values $x$ and $x + \sqrt{x}$. We shall prove it, based on the strictly proven first probabilistic model and relation (4.2).

Thus, we must prove, that:

$$\sqrt{Li(x) - Li^2(x)/x} < \sqrt{x} . \tag{4.6}$$

It is known, that:

$$Li(x) = x/\ln(x) + Cx/\ln^2(x) = x(1/\ln(x) + C/\ln^2(x)), \tag{4.7}$$

where $1 < C < 2$.

If $x \geq e^2$ and $1 < C < 2$, then:

$$1/\ln(x) + C/\ln^2(x) < 1. \tag{4.8}$$

Therefore, for the value $x > e^2$ from (4.7) and (4.8), it follows that $Li(x) = x(1/\ln(x) + C/\ln^2(x)) < x$. So, for the values $x > e^2$ we have that $\sqrt{Li(x) - Li^2(x)/x} < \sqrt{Li(x)} < \sqrt{x}$. Q.E.D.

## 5. PROBABILISTIC MODELS OF THE DISTRIBUTION OF PRIMES IN THE ARITHMETIC PROGRESSIONS

Let us consider the third probabilistic model. Suppose there is an arithmetic progression $f(i) = ki + l$ where $(k, l) = 1$.

Write the values $f(0) = l, f(1) = k + l, f(2) = 2k + l, ... f(n) = nk + l = x$ on different balls, indistinguishable to the touch, put them in a box and mix.

Choose at random the first ball from the box and assign the value $I_1 = 1$ to the random variable of the indicator of success, if the first ball number is prime or the value $I_1 = 0$ to the random variable, if not. Then return the ball in the box, mix the balls and take out the second ball from the box. Assign the value $I_2 = 1$ to the random variable (the indicator of success), if the second ball number is prime or the value $I_2 = 0$ to the random variable if not. Repeat it again $x/k$ times.

Since the balls are returned to the box and the next ball is selected under the same conditions, the random variables $I_i$ are independent and the probability to choose prime in each test $p_i$ is equal to $p_i = p$.

We have already considered the characteristics of the random variable: $M(I_i) = p$, $D(I_i) = p(1-p)$.

The density of the sequence $g(n)$, as a proportion of the sequence $f(n)$ is determined by the formula:

$$P(g/f, 1, n) = \pi(g, 1, n) / \pi(f, 1, n), \qquad (5.1)$$

where $\pi(g, 1, n), \pi(f, 1, n)$ are the numbers of members of the sequence $g(n)$ and $f(n)$, respectively, in the range of the natural numbers $[1, n)$.

Based on [4], the density $P(g/f,1,n)$ is finite probability measure in the range of the natural numbers $[1,n)$. Therefore, $p = P(g/f,1,n)$ is satisfied. The probability of randomly selected ball number to be prime is equal to:

$$p = k/\varphi(k)(1/\ln(x) + o(1/\ln(x))), \qquad (5.2)$$

where $\varphi(k)$ is the Euler function (based on Dirichlet's theorem about the number of the primes in arithmetic progression and (5.1)).

Let us consider the random variable:

$$I(x,k) = \sum_{i=1}^{x/k}(I_i). \qquad (5.3)$$

From the linearity of the expectation and (5.2) we obtain that:
$$M(I(x,k)) = (xp)/k = x(1/\ln(x) + o(1/\ln(x)))/\varphi(k). \qquad (5.4)$$

Based on the independence of random variables $I_i$, the variance of the random variable $I(x,k)$ is equal to:

$$D(I(x,k)) = xp(1-p)/k = \frac{x}{\varphi(k)}(1/\ln(x) + o(1/\ln(x)))[1 - (\frac{k}{\varphi(k)}(1/\ln(x) + o(1/\ln(x))))].(5.5)$$

The random variable $I(x,k)$ in the probabilistic model can be considered, following [1], as the number of members of the arithmetic progression $f(n) = kn + l$, $(k,l) = 1$, which are the primes.

We consider the third probabilistic model for the case, when $o(1/ln(x))$ is the function $f(x) = \sum_{i=1}^{\infty}(i-1)/\ln^i(x) = Li(x)/x - 1/\ln(x)$. Substitute the function $f(x)$ into the formula for the expectation (5.4) of the random variable $I(x,k)$ and obtain:

$$M_3(I(x,k)) = \frac{1}{\varphi(k)}Li(x). \qquad (5.6)$$

Now we substitute the function $f(x)$ in the formula for the variance (5.5) of a random variable $I(x,k)$ and get:

$$D_3(I(x,k)) = \frac{1}{\varphi(k)}Li(x)[1 - \frac{kLi(x)}{\varphi(k)x}] = \frac{Li(x)}{\varphi(k)} - \frac{kLi^2(x)}{\varphi^2(k)x}. \qquad (5.7)$$

The random variable $I(x,k)$ is the sum of the independent random variables $I_i$ with finite variance and thus, it has the binomial distribution. Based on Moivre-Laplace theorem [3] the limit distribution for $I(x,k)$ is normal, therefore the expression is true for large values $x$:

$$P(|I(x,k) - M(I(x,k))| < C\sqrt{D(I(x,k))}) \approx F(C), \qquad (5.8)$$

where $F(C)$ is the function of the module of the standard normal distribution at point $C$.

Substitute the characteristics of the random variable $I(x,k)$ into the expression (5.8) defined by the formulas (5.6) and (5.7) and obtain:

$$P(|I(x,k) - \frac{1}{\varphi(k)} Li(x)| < C\sqrt{\frac{Li(x)}{\varphi(k)} - \frac{kLi^2(x)}{\varphi^2(k)x}}) \approx F(C). \qquad (5.9)$$

Thus, based on the expression (5.9) we can select the value $C$ such that the probability of the event $|I(x,k) - \frac{1}{\varphi(k)} Li(x)| < C\sqrt{\frac{Li(x)}{\varphi(k)} - \frac{kLi^2(x)}{\varphi^2(k)x}}$ is arbitrarily close to 1.

Note that formulas (5.8) and (5.9) are valid for large values $x/k$.

Let us analyze the third probabilistic model. In this model, as in the first one, the ball is returned to the box after it was chosen. Therefore, there is the possibility to choose the same ball a few times in this model. Such event does not happen in the real situation, when one calculates the number of primes of the arithmetical progression $p_i = k/\varphi(k)\ln(i)$ within the interval of the natural numbers from 1 to $i$. Therefore, this probabilistic model requires clarification.

Consider the fourth probabilistic model that is free from the drawback of the third probabilistic model. We will be based on Cramer's model discussed above. Let make a heuristic assumption, that the probability of a member of the arithmetic progression $kn+l, (k,l) = 1, k+l > 2$ to be a prime (in Kramer's model – the probability to choose from $i$-th urn a white ball) is equal to:

$$p_i = k/\varphi(k)\ln(i). \qquad (5.10)$$

The sequence of primes in the arithmetic progression is specified integer, non-negative, strictly increasing, therefore belongs to the class $C$ Kramer's model. Consider the random variable $I(n) = \sum_{i=1}^{n}(I_i)$, where the random variable $I_i$ takes the value 1 if we get a white ball from the $i$-th urn and the value 0 otherwise. Summarizing Kramer random variable $I(n)$ can be regarded as the number of members of the arithmetic progression $f(t) = kt + l, (k,l) = 1$, which are the primes.

We have already determined the characteristics of the random in the second probabilistic model:

$$M(I(n)) = \sum_{i=1}^{n}(p_i), D(I(n)) = \sum_{i=1}^{n}(p_i)(1-p_i).$$

Based on (5.10) we find the expectation of the random variable $I(n)$:

$$M(I(n)) = k/\varphi(k)\sum_{i=1}^{n}\frac{1}{\ln(ki+l)} \approx k/\varphi(k)\int_{t=1}^{n}\frac{dt}{\ln(kt+l)}. \qquad (5.11)$$

The approximation sign in formula (5.11) should be considered as in formulas (2.22) and (2.23), therefore the difference between the sum and the integral is negligible (for large values $t$).

Change the variables in (5.11) $u = kt + l$ and obtain:

$$M(I(n)) \approx k/\varphi(k)\int_{t=1}^{n}\frac{dt}{\ln(kt+l)} = k/k\varphi(k)\int_{u=k+l}^{x}\frac{du}{\ln(u)} = 1/\varphi(k)\int_{u=k+l}^{x}\frac{du}{\ln(u)}. \qquad (5.12)$$

Equation (5.12) can be written as:

$$M(I(n)) \approx 1/\varphi(k)\int_{u=k+l}^{x}\frac{du}{\ln(u)} = 1/\varphi(k)[Li(x) - \int_{2}^{k+l}\frac{du}{\ln(u)}]. \qquad (5.13)$$

The value of the integral in (5.13) is bounded by $\int_{2}^{k+l}\frac{du}{\ln(u)} < (k+l-2)/\ln(2)$, and the following relation is held for large values $x$:

$$M(I(n)) \approx Li(x)/\varphi(k). \qquad (5.14)$$

The variance of the random variable $I(n)$ is equal to:

$$D(I(n)) = k/\varphi(k)\sum_{i=1}^{n}\frac{1}{\ln(ki+l)} - k^2/\varphi^2(k)\sum_{i=1}^{n}\frac{1}{\ln^2(ki+l)} \approx 1/\varphi(k)\int_{k+l}^{x}\frac{du}{\ln(u)} - k/\varphi^2(k)\int_{k+l}^{x}\frac{du}{\ln^2(u)}. \qquad (5.15)$$

The approximation sign in formula (5.15) should be considered in the same sense as in (5.11). As showed earlier, the difference between the sum and the integral is negligible for large values $x$. The value of the integral in (5.15) is bounded by $\int_2^{k+l} \frac{du}{\ln^2(u)} < (k+l-2)/\ln^2(2)$ and the following formula is held for the variance of the random variable $I(n)$ (for large values $x$):

$$D(I(n)) \approx 1/\varphi(k) \int_2^x \frac{du}{\ln(u)} - k/\varphi^2(k) \int_2^{k+l} \frac{du}{\ln^2(u)}. \tag{5.16}$$

The random variable $I(n)$ is the sum of mutually independent random variables $I_i$. Based on Lyapunov theorem the limit distribution for the random variable $I(n)$ is normal, therefore the expression (for large values $n$) is true:

$$P(|I(n) - M(I(n))| < C \cdot D(I(n))) \approx F(C), \tag{5.17}$$

where $F(C)$ is the value of the function of the module of the standard normal distribution at point $C$. Substitute the characteristics obtained for the random variable $I(n)$ (5.14), (5.16) into the expression (5.17) and obtain:

$$P(|I(n) - Li(x)/\varphi(k)| < C\sqrt{Li(x)/\varphi(k) - \frac{k}{\varphi^2(k)} \int_2^x \frac{du}{\ln^2(u)}}) \approx F(C). \tag{5.18}$$

Thus, based on (5.18), you can select such value $C$ that the probability of the event $|I(n) - Li(x)/\varphi(k)| < C\sqrt{Li(x)/\varphi(k) - \frac{k}{\varphi^2(k)} \int_2^x \frac{du}{\ln^2(u)}}$ is arbitrarily close to 1.

### 6. ANALISYS OF THE PROBABILIISTIC MODELS OF THE DISTRIBUTION OF PRIMES IN THE ARITHMETIC PROGRESSION

We have the following dispersion relation for probabilistic models of the distribution of the primes in arithmetic progressions:

$$D(I(n)) < D_3(I(x,k)) < Li(x)/\varphi(k), \tag{6.1}$$

where the value $D_3(I(x,k))$ is defined by (5.7) for the third probabilistic model and $D(I(n))$ is defined by (5.16) for the fourth probabilistic model. Relation (6.1) is proved by the particular case of Cauchy-Schwarz inequality using formula (3.1).

The generalized Riemann conjecture is equivalent to the formula:

$$|\pi(x,k,l) - Li(x)/\varphi(k)| < C_1 \sqrt{x} \ln(x), \tag{6.2}$$

where $\pi(x,k,l)$ is the number of primes not exceeding the number $x$ in the arithmetical progression $kx+l, (k,l)=1$. Comparing (5.9) and (5.18) for the probabilistic models of distribution of primes in the arithmetic progression with formula (6.2) we see that these estimates are stronger. However, as shown earlier, based on Littlewoods's work, these estimates are valid.

On the other hand, formula (6.2) is valid for all $x \geq k$ and formulas (5.9) and (5.18) are valid only with a certain probability. The approximate calculation often use in practice and therefore formulas (5.9), (5.18) are convenient to use in these cases. The value of the constant $C$ in these formulas is selected depending on the desired accuracy.

The second and fourth probabilistic models are conjectures, since they are based on assumptions. The first and the third probabilistic models are the proved assertions.

As mentioned earlier, formulas (5.8) and (5.9) are valid for the value $x/k$ tending to infinity, i.e. when $k = o(x)$. For example, in the case when $k = x^a$ ($0 < a < 1$). This corresponds to Titchmarsh theorem [10]: "Let we have $0 < a < 1$ and $1 \leq k \leq a^x$. Then there exists a constant $C = C(a)$ such that: $\pi(x,k,l) < Cx/\varphi(k)\ln(x)$ for all $(l,k)=1, 0 \leq l < k$" and Theorem [11]: "On the assumption of the previous theorem, for each $k$ we can find such $c_1, c_2$, that there exists more than $c_1\varphi(k)$ different values $l$, for which: $\pi(x,k,l) > c_2 x/\varphi(k)\ln(x)$".

Consider formula (5.9). The number of primes in the arithmetic progression $kt+l, (k,l)=1$ not exceeding $x$ depends on the values $x,k$. The value $Li(x)/\varphi(k)$ is the average value over all $l$ and depends only on the values $x,k$. The value $C\sqrt{D(x,k)}$ depends on the accuracy $(C)$ and values $x,k$.

If we have two arithmetical progressions $kt+l_1, (k,l_1)=1$ and $kt+l_2, (k,l_2)=1$, then from formula (5.9) we get:

$$|I(x,k) - \frac{1}{\varphi(k)}Li(x)| < C\sqrt{\frac{Li(x)}{\varphi(k)} - \frac{kLi^2(x)}{\varphi^2(k)x}} \qquad (6.3)$$

with probability arbitrarily close to 1.

Similarly to the above theorem of Littlewoods [6], it can be proved, that

$$|\pi(x,4,1) - \pi(x,4,3)| < C_4\sqrt{x}\ln\ln\ln(x)/\ln(x). \qquad (6.4)$$

Comparing formulas (6.3) and (6.4) we see that (6.3) is more general than (6.4).

The fourth probabilistic model is more accurate than the third one, as it is free of its drawbacks. On the other hand, based on (6.1), the dispersion of third probabilistic model (5.7) is greater than that of the fourth probabilistic model (5.16) and therefore it can be used to get upper-bound estimate in the inequality (5.18). Since the third probabilistic model is rigorously proved and (5.18) is also rigorously proved assertion.

## 7. APPLICATION OF THE PROBABILISTIC MODELS OF THE DISTRIBUTION OF PRIMES IN THE ARITHMETIC PROGRESSION

Based on the third probabilistic model, there exists at least one prime number belonging to the arithmetic progression $kn+l, (k,l)=1$ between the values:

$$x \text{ and } x+D\sqrt{Li(x)/k - kLi^2(x)/\varphi^2(x)x} \qquad (7.1)$$

for the fixed value $k$, large $x$ and appropriate $D>0$ with the probability close to 1. The condition (7.1) is stronger than the corresponding conditions resulted from the extended Riemann conjecture.

Based on (5.18) we obtain:

$$max_{(k,l)=1} | I(x,k) - Li(x)/\varphi(k) | < C_1\sqrt{Li(x)/\varphi(k) - kLi^2(x)/\varphi^2 x} < C_1\sqrt{Li(x)/\varphi(k)}, \qquad (7.2)$$

where the maximum is taken over all $l$ coprime with $k$.

According to Elliott–Halberstam conjecture for all $0<a<1$ and all $A>0$ there exists such $C>0$ that is held:

$$\sum_{1\leq k\leq x^a} \{max_{(k,l)=1} | \pi(x,k,l) - Li(x)/\varphi(k) |\} \leq Cx/\ln^A(x) \qquad (7.3)$$

for all $x>2$.

This conjecture was proved by Bombieri and Vinogradov for values $a<1/2$.

Assertion 5. Elliott-Halberstam conjecture is held with probability arbitrarily close to 1 for all $0<a<1$.

Proof. Based on the weak estimate (7.2) we need to prove that:

$$C_1\sqrt{Li(x)/\varphi(k)} \leq Cx/\ln^A(x). \qquad (7.4)$$

From $Li(x) = x/\ln(x) + C_2 x/\ln^2(x)$ and (7.4) we see that we must prove the inequality:

$$C_1 \sqrt{x/\ln(x) + C_2 x/\ln^2(x)} \sum_{1 \le k \le x^a} 1/\sqrt{\varphi(k)} \le Cx/\ln^A(x).$$

Based on the Landau theorem [12]:

$$\varphi(n) > C_3 n/\ln(\ln(n)), \tag{7.5}$$

if $n > 2$. Therefore, based on (7.5), we obtain:

$$1/\sqrt{\varphi(n)} < \sqrt{\ln\ln(n)/C_3 n} \tag{7.6}$$

for $n > 2$. Thus, based on (7.6), we obtain:

$$\sum_{1 \le k \le x^a} 1/\sqrt{\varphi(k)} < 2 + 1/\sqrt{C_3} \sum_{2 < k \le x^a} \sqrt{\ln\ln(k)/k} < 2 + 1/\sqrt{C_3} \sqrt{\ln\ln(x^a)} \sum_{2 < k \le x^a} 1/\sqrt{k} \tag{7.7}$$

since $\varphi(1) = \varphi(2) = 1$. Based on assertion 3 we obtain:

$$\sum_{3 \le k \le y} 1/\sqrt{k} = \int_3^y \frac{dt}{t^{1/2}} + C + O(1/\sqrt{y}) = 2\sqrt{y} - 2\sqrt{3} + C + O(1/\sqrt{y}) = O(\sqrt{y}). \tag{7.8}$$

Therefore, the inequality:

$$\sum_{3 \le k \le x^a} 1/\sqrt{k} \le C_4 x^{a/2} \tag{7.9}$$

is satisfied. Substitute the estimate (7.7) in (7.9) and obtain:

$$\sum_{1 \le k \le x^a} 1/\sqrt{\varphi(k)} < 2 + C_4 x^{a/2}/\sqrt{C_3}\sqrt{\ln\ln(x^a)} \sum_{1 \le k \le x^a} 1/\sqrt{\varphi(k)} < 2 + C_4 x^{a/2}/\sqrt{C_3}\sqrt{\ln\ln(x^a)}. \tag{7.10}$$

Therefore, based on (7.10), we obtain the estimate:

$$C_1 \sum_{1 \le k \le x^a} \sqrt{Li(x)/\varphi(k)} < C_1 \sqrt{x/\ln(x) + C_2 x/\ln^2(x)} \cdot (2 + C_4 x^{a/2}/\sqrt{C_3}\sqrt{\ln(a) + \ln\ln(x)}), \tag{7.11}$$

where $0 < a < 1$.

From (7.11) it follows that the estimate is held:

$$C_1 \sum_{1 \le k \le x^a} \sqrt{Li(x)/\varphi(k)} < C \frac{x^{(a+1)/2}(\ln\ln(x))^{1/2}}{\ln(x)^{1/2}} < Cx^{(a+1)/2} \tag{7.12}$$

for values $x > x_1$ where $0 < a < 1$. The following inequality for values $x > x_2$ is held:

$$x^{(1-a)/2} > \ln^A(x) \tag{7.13}$$

for any $A > 0$ and $0 < a < 1$. Therefore, based on (7.12) and (7.13) is held:

$$Cx^{(a+1)/2} < Cx/\ln^A(x) \tag{7.14}$$

for values $x > x_2$ and $C > 0$.

Based on (7.2), (7.12) and (7.14) we obtain the estimate for the values $x > max(x_1, x_2)$:

$$max_{(k,l)=1}\{|I(x,k) - Li(x)/\varphi(k)|\} < C_1 \sum_{1 \leq k \leq x^a} \sqrt{Li(x)/\varphi(k)} < Cx/\ln^A(x). \tag{7.15}$$

Q.E.D.

## 8. CONCLUSIONS AND SUGGESTIONS FOR FURTHER WORK

Elliott-Halberstam conjecture with the probability close to 1 was proved in this paper using the weak estimate (7.2).

However, as shown in [13], from the validity of Elliott-Halberstam conjecture it follows, that there is an infinite number of the prime pairs that differ by no more than 12.

I hope that using the strong estimate (7.2) it will be able to find the probabilistic approach to the solution of the binary Goldbach problem and to prove that there is an infinite number of prime twins.

## 9. ACKNOWLEDGEMENTS

Thank you, everyone, who has contributed to the discussion of this paper.